\documentstyle[10pt,pb-diagram]{article}
\newtheorem{thm}{Theorem}[section]
\newtheorem{de}[thm]{Definition}
\newtheorem{lem}[thm]{Lemma}

\newtheorem{cor}[thm]{Corollary}

\newtheorem{pro}[thm]{Proposition}
\newenvironment{prf}{{\bf Proof}:}{\begin{flushright} $\Box$ \end{flushright}}
\newtheorem{exa}[thm]{Example}

\newcommand{\shoexaseq}[3]{{$ 0 \rightarrow #1 \rightarrow #2 \rightarrow #3 \rightarrow 0$}}
\newcommand{\Shoexaseq}[3]{{$$ 0 \rightarrow #1 \rightarrow #2 \rightarrow #3 \rightarrow 0$$}}
\newcommand{\rar}[1]{\stackrel{#1}{\rightarrow}}

\newcommand{\Ecal}{{\cal E}}
\newcommand{\Fcal}{{\cal F}}
\newcommand{\Kcal}{{\cal K}}
\newcommand{\Ical}{{\cal I}}
\newcommand{\Ncal}{{\cal N}}
\newcommand{\Tcal}{{\cal T}}
\newcommand{\Gcal}{{\cal G}}
\newcommand{\Ocal}{{\cal O}}
\newcommand{\Pcal}{{\cal P}}
\newcommand{\Qcal}{{\cal Q}}
\newcommand{\Rcal}{{\cal R}}
\newcommand{\Hcal}{{\cal H}}

\newcommand{\Id}{{\rm Id}}
\newcommand{\rk}{{\rm rk}}
\newcommand{\Hom}{{\rm Hom}}
\newcommand{\End}{{\rm End}}
\newcommand{\Grass}{{\rm Grass}}

\newcommand{\Pic}{{\rm Pic}}

\newcommand{\HNP}{{\rm HNP}}
\newcommand{\im}{{\rm im}}

\newcommand{\di}{{\rm div}}
\newcommand{\dd}{{\rm d}}
\newcommand{\dual}{^\lor}
\newcommand{\ra}{\rightarrow}
\newcommand{\onto}{{\rightarrow \hspace{-0.8em} \rightarrow}}

\newcommand{\fdop}{{\rm I \hspace{-0.1em} F}}
\newcommand{\pdop}{{\rm I \hspace{-0.1em} P}}

\newcommand{\ndop}{{\rm I \hspace{-0.1em} N}}

\newcommand{\rdop}{{\rm I \hspace{-0.1em} R}}
\newcommand{\cdop}{{\unitlength=1.0pt
\begin{picture}(10.00,5.00)(152.00,677.00)
\put(156.00,683.00){\line(0,-1){6.00}}
\put(156.00,680.00){\makebox(0,0)[cc]{\bf{C}}}
\end{picture}}}

\begin{document}

\title{THE SPECTRUM OF SEMISTABLE VECTOR BUNDLES ON CERTAIN ALGEBRAIC SPACES}
\author{H. Kurke \and O. Teschke}
\maketitle

\

\tableofcontents
\setcounter{section}{-1}
\section{Introduction}

For the variety $Z=\pdop^3$ the spectrum of semistable sheaves was studied by several authors, starting with \cite{B-E} and
\cite{Ba}. The aim was to control - to a certain extend - the first cohomology of the sheaf twisted
with $\Ocal(l)$ ($l=-1,-2,\dots$). Ein and Sols showed the existence of a spectrum for reflexive sheaves of rank two on 
quadrics (\cite{E-S}) and proved the connectedness of 
the spectrum. 
\\
\\We have in mind to adapt these methods to a class of compact algebraic spaces, where algebraic twistor spaces are the main examples.
Combining results of Poon (\cite{Po}) and Campana (\cite{Ca}) such twistor spaces are twistor spaces of self-dual conformal manifolds
homeomorphic to $S^4$ or connected sums of $n$ complex projective spaces. For each $n$ there are explicit families of such algebraic twistor spaces (\cite{LeB}, \cite{Ku}).
\\In this paper, we actually do not use too much of the property of being a twistor space. The only properties we need are the vanishing
of the first cohomology of the structure sheaf, the existence of a half-anti-canonical sheaf, and certain properties of the Chow ring.
\\
\\Except of the spaces $\pdop^3$ (twistor space of $S^4$) and the flag variety (twistor space of $\pdop^2_{\cdop}$), this spaces are not
K\"ahlerian. Thus we need an extended notion of stability in this context (section 1). 
\\The usual notion of "stability of torsion-free sheaves" on a projective variety takes into account properties of subsheaves restricted
to the family of curves which occur as linear sections of a fixed projective embedding.
\\Instead of this family we consider an arbitrary family of smooth curves on our space $Z$, parametrized by a smooth connected algebraic
space $H$, subject to the following conditions:
\\(i)  The curves passing through an arbitrary given point of $Z$ are parametrized by a smooth subspace of positive dimension, 
independent on the given point;
\\(ii) generically, this subspace is connected.
\\Using this properties, we can show that the analogue of the theorem of Grauert-M\"ulich holds in our situation (section 2).
\\
\\For complete families of curves (in the sense of Kodaira-Spencer) the first property is satisfied if through each point $P$ of $Z$ passes 
one curve and its normal bundle $\Ncal$ satisfies $H^1(\Ncal(-P))=0$. The second property required for our family seems to be quite delicate and not so easy to verify (see section 3). 
\\ 
\\The second point where the projective embedding was used originally was the study of the twisted cohomology (with the twisting sheaf
$\Ocal(1)$).  In our case the twisting sheaf is replaced by the 
half-anticanonical line bundle. If we specialize this to the known case of $\pdop^3$, it turns out that the spectrum in our sense 
controls only the differences $h^1(\Fcal(-1))-h^1(\Fcal(-3))$, $h^1(\Fcal(-3))-h^1(\Fcal(-5)), \dots$. Comparing the known results in 
this case with our results, one could hope to prove stronger properties about the connectedness of the spectrum than we obtained so far.
\\
\\{\bf Acknowledgements:} Parts of this work were done by the second author during PRAGMATIC 1998 in Catania. He likes to thank Europroj, 
the $Graduiertenkolleg$ "Geometry and Non-Linear Analysis" in Berlin and the University of Bayreuth for financial support.
Furthermore, he likes to thank the organizers of PRAGMATIC 1998  for the 
warm atmosphere during this school, and the lecturers Lawrence Ein and Vladimir Masek for useful discussions and valuable hints. 
The valuable hints from this universit
\section{Relative stability}

From now on, we assume $Z$ to be a smooth connected compact complex algebraic space, $H$ to be a smooth connected 
complex algebraic space, and $W\subset Z\times H$ a family of smooth curves on Z. 
\\Let $\mu , \nu$ be the projections 
$$\begin{diagram}
\node{} \node{W\subset Z\times H} \arrow{sw,l}{\mu} \arrow{se,t}{\nu}
\\ \node{H} \node{} \node{Z} \node{.}
\end{diagram}$$
\\Moreover, we will assume 
$$(i)\quad W\rar{\nu}Z \mbox{\quad is smooth, and }\qquad \qquad$$
$$(ii) \quad\mbox{the generic fibre of } \nu \mbox{ is connected.}$$
\\The second condition is quite technical and often difficult to prove, but we see no possibility to drop it.
\begin{de} Let $\Fcal $ be a torsion free coherent sheaf on $Z$. We define the W-relative
degree and the W-relative slope as 
$$\deg (\Fcal)=\deg (\Fcal |W_h) \mbox{\qquad for $h\in H$ generic}$$ 
and
$$\mu(\Fcal)=\frac{\deg \Fcal}{\rk \Fcal}. \qquad \qquad \qquad \qquad\qquad \qquad$$
$\Fcal$ is $W$-relatively stable (semistable), if for every torsion free coherent subsheaf ${\cal F}'$ 
with $0<rk{\cal F}'<{\cal F}$ there holds: $\mu({\cal F}')<\mu ({\cal F})$  
(resp. $\mu({\cal F}')\le\mu ({\cal F})$).
 \end{de}
Since the family $W$ will be clear from the context from now on, we will omit the word "$W$-relative" in the following.

\begin{pro} If $\Fcal$ is stable, then $\Fcal$ is simple. \end{pro}
\begin{prf} Given $\phi\in \End(\Fcal)$, we consider $\Kcal=\ker(\phi)$,  $\Ical=\im (\phi)$. If we assume 
$\Kcal\neq 0$ and $\phi\neq 0$, then $\mu (\Kcal)<\mu(\Fcal)$, $\mu (\Ical)<\mu(\Fcal)$, hence
$\deg(\Kcal)+\deg{\Ical}<\deg{\Fcal}$, which is impossible.
\\Therefore for all $\lambda\in\cdop$, $\phi - \lambda \Id$ is injective or the zero map. 
\\The first case 
implies that $\det (\phi - \lambda Id):\det\Fcal\ra\det\Fcal$ is an isomorphism. But we can choose $x$ such that $\Fcal _x$ is free
and $\lambda$ such that it is an eigenvalue of $\Fcal(x)\ra \Fcal(x)$. Then $\det(\phi - \lambda \Id)$ is not an isomorphism in $x$,
hence $\phi=\lambda \Id$. \end{prf}
As in \cite{Sha} we can now introduce the Harder-Narasimhan filtration: 
\begin{thm} A torsion free coherent sheaf has a filtration 
$$\Fcal_0=0\subset \Fcal_1\subset \dots \subset \Fcal_m=\Fcal$$
uniquely determined by the properties
\\(i) $\Fcal_{j+1}/\Fcal_j$ is torsion free and semistable of slope $\mu_{j+1}$, $j=0,\dots ,m-1$
\\(ii) $\mu_1>\mu_2>\dots >\mu_m$.
 \end{thm}
\begin{prf}\cite{Sha} \end{prf}
As a corollary from this theorem, there are only finitely many points in $\rdop^2$ in the set
$$B_{\Fcal}=\{(\rk(\Fcal '),\deg(\Fcal '))|\Fcal ' \subset \Fcal\}$$ lying above the line through $(0,0)$ with slope $\mu_{\Fcal}$.
Hence, the convex hull of $B_{\Fcal}$ has an upper boundary consisting of a finite number of edges of slope $\mu_j(\Fcal)$
(with the vertices corresponding to the sheaves in the filtration). This boundary polygon is called the Harder-Narasimhan polygon
(HNP). For families of torsion free sheaves, it is constant on an open subset:

\begin{thm}If $(\Fcal_s)_{s\in S}$ is an algebraic family of torsion free sheaves, there exists a polygone $P_0$ 
and a dense open subset $U\subset S$ such that
\\(i)  $ \HNP(\Fcal_s)=P_0$  for  $s\in U $
\\(ii) $\Fcal$ has a filtration by coherent subsheaves which is the Harder-Narasimhan filtration over $U$.
\\Moreover: For any convex polygone P with vertices in $\ndop^2$ the set  
$$\{s\in S|\HNP(\Fcal_s)\ge P\}$$
is closed. \end{thm}
\begin{prf}\cite{Sha} \end{prf}
We will define now the most important invariant of $W$ in this context, which controls the behavior of
the HNP of semistable sheaves after the restriction to curves in $W$:
\begin{de} $$\dd(W)=\mu_{max}((\Omega^1_{W/Z})_h) \mbox{ for generic } h\in H.$$
\end{de}
However, for concrete examples this definition is quite difficult to handle. But we can do the computation of $\dd (W)$
in terms of the normal bundle of the curves under a small additional assumption:
\begin{pro}If the family W is complete in $h$ for generic $h\in H$, then 
$$\dd (W)=-\mu_{min}(\ker(\Ocal_{W_h}\otimes H^0(\Ncal_{W_h/Z})\ra \Ncal_{W_h/Z}))  $$\end{pro}
\begin{prf} For each $h\in H$ we have the characteristic map 
$$\rho_h:\Tcal_h(H)\ra H^0(\Ncal_{W_h/Z})$$
and $W$ is called complete, if $\rho_h$ is surjective. It induces a map 
$$\Ocal_{W_h}\otimes \Tcal_h(H)\ra \Ocal_{W_h}\otimes H^0(\Ncal_{W_h/Z})\ra \Ncal_{W_h/Z}$$
and if it is complete then 
$$\mu_{min}(\Ocal_{W_h}\otimes H^0(\Ncal_{W_h/Z}))=\mu_{min}(\Ocal_{W_h}\otimes \Tcal_h(H)\rar{1\otimes\rho}\Ncal_{W_h/Z}).$$
If we consider $W\subset Z\times H$ and the corresponding standard exact sequence 
$$0\ra \Ncal^{\dual}_{W/Z\times H}\ra \Omega^1_{Z\times H/H} \otimes\Ocal_W\ra\Omega^1_{W/Z}\ra 0,$$
dualize and restrict to the fibre $W_h$, we get the exact sequence 
$$0\ra(\Theta_{W/Z})|_{W_h}\ra \Ocal_{W_h}\otimes\Tcal_h(H)\rar{1\otimes \rho}\Ncal_{W_h/Z},$$
which proves the claim. \end{prf}

\section{A relative Grauert-M\"ulich theorem}
Using the invariant defined in the previous section, we can now state a Grauert-M\"ulich theorem for families of curves.
\begin{thm}[Grauert-M\"ulich] If $(Z,W,H)$ satisfies the conditions (i),(ii) and $\Fcal$ is a semistable torsion free sheaf 
on $Z$, then for generic $h\in H$ 
$$\mu_j(\Fcal_h)-\mu_{j+1}(\Fcal_h)\le\dd(W)$$ \end{thm}
\begin{prf} Consider the relative HNF of $\nu^*\Fcal$: 
$$0=\tilde{\Fcal}_0 \subset \tilde{\Fcal}_1\subset\dots \subset\tilde{\Fcal}_m=\nu^*\Fcal.$$
According to theorem 1.4. there is a dense open $U\subset H$ such that $\nu^*\Fcal|\mu^{-1}U$ is locally free and the filtration
induces the HNF on $(\nu^*\Fcal)_h=\Fcal|_{W_h}$.
\\Assuming $\mu_j(\Fcal_h)-\mu_{j+1}(\Fcal_h)>\dd(W)$ for some $j$ we will construct a subsheaf $\Gcal$ which violates the semistability
condition. 
\\Consider $X=\mu^{-1}U\subset Z\times U$ and $Y=\nu(X)\subset Z$ (open subspace). Then the map $\nu :X\ra Y$ is smooth and the generic fibre 
of $\nu$ is connected. 
\\We need now the following lemma:
\begin{lem}Let $X\rar{\nu} Y$ be a smooth morphism of complex manifolds with connected generic fibre, $\Fcal$ a locally free sheaf on $Y$ 
and $\tilde{\Gcal}\subset\nu^*\Fcal$ a locally free subsheaf. 
\\If $$\Hom(\Theta_{X/Y},{\cal H}{\rm om}(\tilde{\Gcal},\nu^*\Fcal/\tilde{\Gcal}))=0,$$
then there exists a coherent subsheaf $\Gcal\subset\Fcal$ such that $\tilde{\Gcal}=\nu^*\Gcal$. \end{lem}
\begin{prf} If $p=\rk(\nu^*\Fcal/\tilde{\Gcal})$, then $\tilde{\Gcal}$ defines a morphism $X\rar{f} G=\Grass_p(\Fcal)$ (such that
$\tilde{\Gcal}$, $\nu^*\Fcal/\tilde{\Gcal}$ are induced from the universal subbundle and the quotient bundle, respectively).
\\Therefore $$T(f):\Theta_{X/Y}\ra f^*\Theta_{G/Y}={\cal H}{\rm om}(\tilde{\Gcal},\nu^*\Fcal/\tilde{\Gcal})=0, \mbox{\quad hence
\hspace{0.1cm}} T(f)=0.$$
Hence, there exists a section $s:Y\ra G$ with $f=s\circ\nu$ possibly except from a closed subset of codimension two. 
But a birational morphism of smooth varieties and an exceptional set of codimension at least two is an isomorphism, hence
 $f=s\circ\nu$ holds everywhere. 
\\Therefore, we have $\tilde{\Gcal}=\nu^*\Gcal$, $\Gcal=s^*\tilde{\Gcal}$.
 \end{prf}
Now, we apply this lemma to $\tilde{\Gcal}=\Fcal_j$:
\\If there is a non-zero homomorphism $\alpha: \Theta_{W/Z}\ra{\cal H}{\rm om}(\Fcal_j,\nu^*\Fcal/\Fcal_j)$, then there exist
subquotients $\Pcal,\Qcal,\Rcal$ of the relative HNF of $\Theta_{W/Z}$, $\Fcal_j$ and $\nu^*\Fcal/\Fcal_j$ respectively, such that
$\alpha$ induces a non-zero homomorphism $\Pcal\rar{\alpha}{\cal H}{\rm om}(\Qcal,\Rcal)$.
\\Since ${\cal H}{\rm om}(\Qcal,\Rcal)$ is semistable of slope $\mu(\Rcal)-\mu(\Qcal)=\mu_k-\mu_i$ $(k>j,i\le j)$ it follows
$$\mu_{min}(\Theta_{W/Z})\le\mu(\Pcal)\le \mu_k-\mu_i\le \mu_{j+1}-\mu_j,$$ hence
$$\mu_j-\mu_{j+1}\le -\mu_{min}(\Theta_{W/Z})=\dd(W).$$
Therefore $Hom((\Theta_{W/Z}),{\cal H}{\rm om}(\Fcal_j,\nu^*\Fcal/\Fcal_j))=0$ and $\Fcal_j$ is induced from a coherent subsheaf 
$\Gcal\subset\Fcal$(at least on $Y$). 
\\We can extend $\Gcal$ to a coherent subsheaf of $\Fcal$, then $\Gcal$ violates the semistability condition.
\end{prf}

\section{Examples}
In this section we will consider four examples of families of curves, verify the conditions for the application of the 
Grauert-M\"ulich theorem 2.1 and compute the invariant $\dd(W)$. The common property of the examples $3.1$,$3.2$ and $3.4$ is that
they admit the possibility of a geometric construction of a spectrum for bundles of rank two in a natural way. For the 
first example this was done in\cite{B-E} and for the second in \cite{E-S}, both for families of lines. The construction for the
special families of elliptic curves from $3.4$ is essentially new and will be done in the next section.
\begin{exa}[Lines on $\pdop^n$] {\rm Let $W\subset\pdop^n\times\Grass(1,n)$ be the family of lines on $Z=\pdop^n$. 
The normal bundle of a line
L is $\Ocal_{L}(1)^{n-1}$, hence $$\dd=-\mu_{min}(\ker(\Ocal_L^{2n-2}\ra\Ocal_L(1)^{n-1}))=-\mu_{min}(\Ocal_L(-1)^{n-1})=1,$$
and  the fibres of $\nu$ (the lines through a fixed point) are isomorphic to $\pdop^{n-1}$, hence connected.
The application of the theorem 2.1 gives the classical Grauert-M\"ulich theorem: 
Semistable torsion free sheaves on $\pdop^n$ restricted to a generic line are isomorphic to 
$\bigoplus_{j=m}^n \Ocal_{\pdop^1}(j)^{m_j}$
such that $m_j\neq 0$ for all $m\le j\le n$.}
\end{exa}

\begin{exa}[Lines on quadrics] {\rm Let $Z\subset \pdop^{n+1}$ be a quadric with $n=\dim Z\ge3$, and 
let $W\subset Z\times H$ be the family of lines. The normal bundle is $\Ncal_{L/Z}=\Ocal_L\oplus\Ocal_L(1)^{n-2}$, hence 
$$\dd=-\mu_{min}(\ker(\Ocal_L\otimes H^0(\Ncal_{L/Z})\ra\Ncal_{L/Z}))=-\mu_{min}(\Ocal_L(-1)^{n-2})=1.$$
The connectedness of the fibres of $\nu$ can be seen as follows: 
We can regard $Z$ as a double cover of $\pdop^n$ branched over a quadric $Q\subset\pdop^n$. The lines in $Z$ correspond to
lines in $\pdop^n$ tangent to $Q$, which form a $2(n-1)-2=2n-3$-dimensional family. The fibres of $\nu$ correspond to the family of
tangent lines to $Q$ through a fixed point, which form a connected family for $n\ge 3$.
\\Again, theorem 2.1 ensures that semistable torsion free sheaves on $Z$ restricted to a generic line are isomorphic to 
$\bigoplus_{j=m}^n \Ocal_{\pdop^1}(j)^{m_j}$
such that $m_j\neq 0$ for all $m\le j\le n$.
}
\end{exa}

\begin{exa}[Twistor lines on twistor spaces] {\rm We consider now twistor spaces over self-dual manifolds of positive type homeomorphic to
$S^4$ or $\#_n \pdop^2$ ($n\le 3$). As we have mentioned in the introduction, the crucial point is to check the property (ii) (i.e.
generic connectedness) for the family of the twistor lines. 
\\This is easy for $\pdop^3$ and $\fdop(1,2)$. But for $n=2,3$ the twistor spaces are non-projective, even not K\"ahlerian. We give here
the idea for the more complicated case $n=3$:
\\Let $Z$ be a twistor space over the connected sum of three $\pdop^2$. First, we assume that $|\omega^{-1/2}|$ has base points.
\\In this case, it was shown in \cite{LeB} that $Z$ is a small resolution of a conic bundle over $\pdop^1\times \pdop^1$. 
We will use here the following description of the twistor lines from \cite{Ku} :
$Z$ is a small resolution $Z\ra Z'$ of  
$$Z'\subset \pdop(\Ecal)\ra Q=\pdop^1\times \pdop^1, \quad \Ecal=\Ocal_Q z_0\oplus \Ocal(-2,-1) z_1\oplus\Ocal(-1,-2) z_2$$
given by the equation $$z_1 z_2-\phi_1\phi_2\phi_3 z_0^2=0 \mbox{\quad with \quad} \phi_1 ,\phi_2 ,\phi_3 \mbox{\quad bilinear forms}.$$
Furthermore, we define the surfaces $S_j=\{ z_1=\phi_j=0\}$ and $\bar{S}_j=\{ z_2=\phi_j=0\}$  $(j=1,2,3)$. Then, the twistor lines on $Z$ 
correspond to lines $F$ in $Z'$ with the properties $F.S_j=1,F.\bar{S}_j=1$ $(j=1,2,3)$, which form a smooth connected family $W$ 
(\cite{KK}, \cite{Ku}). 
\\The idea to show the connectedness of a generic fibre of $\nu :W\ra Z$ is the following:  A generic $F$ is given by 
$$y=Sx,\quad z_1=u(x)z_0,\quad z_2=v(x) z_0$$
with $(x,y)=(x_0:x_1,y_0:y_1)$ are coordinates on $\pdop^1\times\pdop^1$, $S\in Sl_2(\cdop)$, and $u(x), v(x)$ are cubic forms
with $$u(x)v(x)=\phi_1(x,Sx)\phi_2(x,Sx)\phi_3(x,Sx)$$ and the property, that $u(x)$ (hence also $v(x)$) has a common factor with each $\phi_j$
(this comes from $F.S_j=1$). Hence, we can describe the fibres of the morphism (i.e. the twistor lines through a fixed point ($\bar{x},\bar{y}$))
as a surface $$V \subset \pdop^2(q_0:q_1:q_2)\times \pdop^3(c_0:c_1:c_2:c_3)$$ (with $q_j$ being the coefficients of quadric und $c_j$ the coefficients
of cubic forms). More precisely, $V$ is given as the intersection of three hypersurfaces defined by resultants. By sending a point in $V$ 
to the
common factor of both the quadric and the cubic form, we obtain a fibration $V\ra \pdop^1$, with reducible but generically connected fibres
(this is quite technical to show, so we omit the proof here). Hence $V$ is connected.
\\Let's regard now the second case, i.e. assume $|\omega^{-1/2}|$ is base-point free. Then $Z$ is a small resolution of a double solid
branched over a 13-nodal quartic (\cite{KK}). It was shown in \cite{Ca} that every twistor space of the first type (i.e. the "conic bundle
type") occurs as a deformation of a twistor space of the second type (i.e. the "double solid type"). Hence, by the semicontinuity 
theorem, we have also the connectedness of the generic fibre in the second case. 
\\Since the normal bundle of twistor lines is $\Ocal_F(1)^2$, we obtain again $\dd (W)=1$ and the familiar type of a Grauert-M\"ulich theorem.
}\end{exa}

\begin{exa}[Families of elliptic curves]{\rm While the curves in the previous examples were rational, we consider now families of elliptic curves for a quite 
large class of algebraic spaces. 
\\Let $Z$ be a threefold with $H^1(K_Z)=0$. Assume the existence of a root $L\in\Pic (Z)$ of the anticanonical bundle $L^{\otimes 2}
=K_Z^{-1}$ with the following properties:
\\(i) \hspace*{0.12cm}  $H^1(L^{-1})=0$
\\(ii)\hspace*{0.1cm}  $|L|$ is free, dim$|L|\ge 3$
\\(iii) $H\subset \Grass(2,H^0(L))$ is an open subset such that the base curve of the corresponding pencil is smooth.
\\ \\So $H\subset \Grass(2,H^0(L))$ parametrizes a family $W\subset Z\times H$ of base curves. Under the previous assumptions the 
following holds: 
\begin{pro}  $W\subset Z\times H$ is a generically complete family of elliptic curves with 
$$\dd (W)=\frac{\dim |L|-1}{\dim |L|-2}$$
\end{pro} 
\begin{prf} Let $\Ical\subset\Ocal_Z$ be the ideal sheaf of $C$, then the sequence
\Shoexaseq { L^{-\otimes 2} }{ L^{-1}\oplus L^{-1 }}{ \Ical }
is exact. Using the long exact sequence, Serre duality and (i), we have immediately
$$H^1(\Ical)\simeq H^2(L^{-\otimes 2})\simeq H^2(K_Z)\simeq H^1(\Ocal)\dual =0,$$ hence $H^0(\Ocal_C)=\cdop$, $C$ is connected.
\\$K_C=K_Z|C\otimes \det \Ncal_{C/Z}=\Ocal_C$ by the adjunction formula, hence $C$ is elliptic.
\\For the proof of the completeness we compute the characteristic map:
\\Let $\phi_0,\ ,\phi_1,\psi_1,\dots ,\psi_n$ be a basis of $H^0(L)$ such that $C=\di (\phi_0).\di(\phi_1)$, then $H\subset\cdop^{2n}$
is locally given by 
$$\Phi_0=\phi_0+x_1 \psi_1+\dots +x_n\psi_n \quad \mbox{and} \quad \Phi_1=\phi_1+y_1 \psi_1+\dots +y_n\psi_n$$
Then $\rho : \cdop^{2n} \ra H^0(\Ncal_{C/Z})=H^0(L_C^{\oplus 2})$ is defined by
$$\rho(\frac{\delta}{\delta x_j})={\psi_j|C \choose 0} \quad , \quad \rho(\frac{\delta}{\delta y_j})={0 \choose \psi_j|C},$$
hence $\rho$ is an isomorphism.
\\For the computation of d$(W)$ we need the following lemma:
\begin{lem}Let $C$ be an elliptic curve, $L$ a line bundle of degree $n\ge 2$ and $\Rcal=\ker(\Ocal_C\otimes H^0(L)\ra L)$. 
Then $\Rcal$ is stable, $\mu(\Rcal)=-\frac{n}{n-1}$.
\end{lem}
\begin{prf} We can write $L=M^{\otimes n}, \deg M=1$, and we proceed by induction on n, using 
\Shoexaseq{\Rcal_2=M^{-2}}{\Ocal_C\otimes H^0(M^{\otimes 2})}{M^{\otimes 2}} 
for $n=2$ and 
$$\begin{diagram}
\node{} \node{0} \arrow{s}  \node{0} \arrow{s}  \node{0} \arrow{s} \node{}
\\ \node{0} \arrow{e} \node{\Rcal_n} \arrow{s} \arrow{e} \node{\Ocal_C\otimes H^0(M^{\otimes n})} \arrow{s} \arrow{e}
 \node{M^{\otimes n}} \arrow{s}  \arrow{e}\node{0}
\\ \node{0} \arrow{e} \node{\Rcal_{n+1}} \arrow{s}  \arrow{e}\node{\Ocal_C\otimes H^0(M^{\otimes n+1})} \arrow{s} \arrow{e}
 \node{M^{\otimes n+1}} \arrow{s}  \arrow{e}\node{0}
\\ \node{0} \arrow{e} \node{M^{-1}} \arrow{s}  \arrow{e}\node{\Ocal} \arrow{s}  \arrow{e}\node{\cdop(p)} \arrow{s} \arrow{e}\node{0}
\\ \node{} \node{0} \node{0}  \node{0} \node{}
\end{diagram}$$
with  $M=\Ocal(P)$, for the induction step. We need again a technical sublemma:
\begin{lem}
Let \shoexaseq{\Rcal '}{\Rcal}{M^{-1}} be a nonsplitting sequence with $\deg{\Rcal}=(\rk(\Rcal)+1)$, $\deg M=1$.
\\Then stability of $\Rcal '$ implies stability of $\Rcal$.
\end{lem}
\begin{prf} We have $\mu (\Rcal)=-1-\frac{1}{\rk(\Rcal)}$ and  $\mu (\Rcal ')=-1-\frac{1}{\rk(\Rcal)-1}$.
\\Suppose $\Rcal$ is not stable, then there exists a stable quotient sheaf $\Rcal \onto \Qcal$ with $\rk{\Qcal}\le\rk{\Rcal}$ and 
$\mu{\Qcal}\le\mu{\Rcal}$. Restriction to $\Rcal '$ gives a diagram 
$$\begin{diagram}
\node{\Rcal '} \arrow{s} \node{\subset} \node{\Rcal}\arrow{s}
\\ \node{\Qcal '} \node{\subset} \node{\Qcal}
\end{diagram}$$
with $\mu(\Rcal ')\le \mu(\Qcal ')\le \mu({\Qcal})\le \mu({\Rcal})$. Since $\rk (\Qcal)\le\rk (\Rcal)-1=\rk (\Rcal ')$, which is only
possible if $\Rcal '\simeq\Qcal$, hence the sequence would split.
\end{prf}
We have $H^0(M\otimes\Rcal_{n+1})=\ker(H^0(M)\otimes H^0(M^{\otimes n+1})\ra H^0(M^{\otimes n+2}))=0$, hence the
sequence \shoexaseq{\Rcal_n}{\Rcal_{n+1}}{M^{-1}} does not split. So we can apply the sublemma and proceed by induction.
\end{prf}
From the Riemann-Roch theorem we get $$\deg L_C=h^0(L_C)=h^0(L)-2=\dim |L|-1\ge 2.$$ The application of the lemma finishes the proof.
\end{prf}
There is a lot of examples for $(Z,L)$ with the required properties, for instance ($\pdop^3,\Ocal(2)$) and 
Fano threefolds of index 2 like the cubic in $\pdop^3$,
the Flag manifold $\fdop(1,2)$, the intersection of two quadrics in $\pdop^5$ and the double covering of $\pdop^3$, ramified 
along a quartic surface. On the other hand, also the twistor spaces mentioned in the previous example fit into this setup. Note, however,
that the family of elliptic curves we get from $\omega^{-1/2}$ is different from the family of twistor lines in example $3.3$.
\\The application of theorem $2.1$ to this class of spaces gives the following Grauert-M\"ulich type corollary:
\begin{cor} Let $\Ecal$ be a semistable vector bundle of rank two on a threefold $Z$ (with the previous assumptions). Then 
for a generic elliptic curve $C\in H$, either $\Ecal |C$ is semistable, or $\Ecal |C=L_1\oplus\L_2$ with $\deg L_1\ge\deg L_2\ge\deg L_1-1$.
\end{cor}
{\begin{flushright} $\Box $ \end{flushright}}
} \end{exa}

\section{Spectra of rank two bundles}
\subsection{Construction of the spectrum}
We will give here the construction for the spectrum of vector bundles of rank two on the threefolds from the fourth example of the
previous section. The construction is similar to the construction for lines on $\pdop^1$ given in \cite{B-E}.
\begin{thm} Let $\Ecal$ be a semistable vector bundle of rank 2 on $Z$, normalized such that 
$$1-\mu(L)\le\mu(\Ecal)<0.$$
There exists a locally free sheaf $\Hcal$ on $\pdop^1$ such that
$$(i) \hspace*{0.6cm} H^1(Z,\Ecal (l))\simeq H^0(\Hcal(l)) \mbox{\quad for \quad} l\le0 \quad $$
$$(ii)\hspace*{0.5cm} H^2(Z,\Ecal (l))\simeq H^1(\Hcal(l)) \mbox{\quad for \quad} l\ge -1 $$ 
where $\Ecal(m)=\Ecal\otimes L^{\otimes m}$. 
\\Furthermore,
\begin{eqnarray}r = \rk (\Hcal) & = & -(\mu(\Ecal)+\frac{p_1(\Ecal)c_1(L)}{2}+2) \nonumber \\
 d = \deg(\Hcal) & = & \frac{3 c_1(\Ecal)c_2(\Ecal)-c_1(\Ecal)^3}{6}-\frac{c_1(\Ecal)c_2(Z)}{12}+\frac{\mu(\Ecal)}{3}\nonumber 
\end{eqnarray}
 \end{thm}
\begin{prf}  
We choose a generic pencil $\Lambda\subset |L|$ with base curve $C$. Denoting the blowing-up of $Z$ along $C$ as 
$\hat Z=Bl_C Z$, we have the following morphisms
$$\begin{diagram}
\node{} \node{\hat Z} \arrow{sw,l}{\beta} \arrow{se,t}{f}
\\ \node{Z} \node{} \node{\pdop^1}
\end{diagram}$$
with $f$ induced by $\Lambda$. 
\\Then, we define $\Hcal=R^1 f_* (\beta^*\Ecal)$.
\\The fibres of $f$ are the surfaces $S\in\Lambda$, and we have exact sequences \Shoexaseq{\Ecal_S(-l-1)}{\Ecal_S(-l)}{\Ecal_C(-l)}
Since $\Ecal_C$ has no sections, we get $H^0(\Ecal_S(-l-1))\simeq H^0(\Ecal(-l))$ for $l \ge 0$, hence $f_*(\beta^*\Ecal(-l))=0$ for 
$l \ge 0$.
\\By Serre duality, we have also $H^2(\Ecal_S(l))=H^0(\Ecal_S\dual(-l-1))\dual=0$ for $ l\ge 0$, hence $R^2 f_*(\beta^*\Ecal(-l))=0$
for $l\ge 0$.
\\Therefore, $R^1 f_* (\beta^*\Ecal)=\Hcal$ is a locally free sheaf compatible with base change. 
\\If we denote $\lambda=c_1(L)$, hence $2\lambda =c_1(Z)$, we have 
\begin{eqnarray}
r  =  \rk\Hcal   & = & h^1(\Ecal_S) = -\chi(\Ecal_S) \nonumber
\\               & = & -(-c_2(\Ecal)\lambda + \frac{1}{2}c_1(\Ecal)^2 \lambda +\frac{1}{2} c_1(\Ecal) \lambda^2 + 2)  \nonumber
\\               & = & -\frac{p_1(\Ecal)\lambda}{2} - \mu{\Ecal} - 2 \nonumber
\\d = \deg(\Hcal)& = & \chi(\Hcal) - \rk(\Hcal) = - \chi(\Ecal) - r  \nonumber
\end{eqnarray}
Using 
\begin{eqnarray}
Todd_1(Z) & = & \lambda \nonumber\\
Todd_2(Z) & = & \frac{c_1(Z)^2+c_2(Z)}{12} = \frac{\lambda^2}{3}+\frac{c_2(Z)}{12} \nonumber \\
Todd_3(Z) & = & \frac{c_1(Z)c_2(Z)}{24}    = \frac{\lambda c_2(Z)}{12} = 1, \nonumber
\end{eqnarray}
we apply now the Riemann-Roch theorem on $Z$ to compute
\begin{eqnarray}
\chi(\Ecal) & = & \frac{c_1(\Ecal)^3-3 c_1(\Ecal^)3 c_2(\Ecal)}{6} + \frac{p_1(\Ecal)}{2}\lambda + 
c_1(\Ecal) (\frac{\lambda^2}{3}+\frac{c_2(Z)}{12}) +2 \nonumber\\
            & = & \frac{c_1(\Ecal)^3-3 c_1(\Ecal^)3 c_2(\Ecal)}{6} + \frac{p_1(\Ecal)}{2}\lambda + 
\frac{2}{3}\mu(\Ecal) + \frac{c_1(\Ecal)c_2(\Ecal)}{12} + 2, \nonumber 
\end{eqnarray}
hence
$$ d = \deg(\Hcal) = \frac{3 c_1(\Ecal^)3 c_2(\Ecal)-c_1(\Ecal)^3}{6} - \frac{c_1(\Ecal)c_2(\Ecal)}{12} + \frac{\mu(\Ecal)}{3} $$
We finish the proof with the computation of the $H^1(Z,\Ecal(l))$:
\\Since we have the vanishing of $H^0(\Ecal_S(l))$ or $H^2(\Ecal_S(l))$ for $l\le 0$ or $l\ge -1$ respectively, the spectral sequence 
$$ H^i({\pdop}^1, R^j f_*(\beta^* \Ecal(l)))\Rightarrow H^{i+j}(\hat{Z},\beta^* \Ecal(l))\simeq  H^{i+j}(Z,\Ecal(l))$$
degenerates, which gives us isomorphisms 
$$ H^0({\pdop}^1, R^1 f_*(\beta^* \Ecal(l))) \simeq H^1(Z,\Ecal(l)) \mbox{\quad for \quad} l\le 0 \quad$$
$$ H^1({\pdop}^1, R^2 f_*(\beta^* \Ecal(l))) \simeq H^2(Z,\Ecal(l)) \mbox{\quad for \quad} l\ge -1.$$
So it remains to show that
$$ H^{0;1}({\pdop}^1, R^1 f_*(\beta^* \Ecal(l)))\simeq  H^{0;1}({\pdop}^1, R^1 f_*(\beta^* \Ecal)(l))=H^{0;1}(\Hcal(l)).$$
We can restrict the proof to the $H^0$ case, the other one is analogous.
Let $E\subset \hat Z$ be the exceptional locus of the blow-up. Then of $\beta^*l\simeq f^*\Ocal (1)(E)$, hence
$$R^1 f_* (\beta^*\Ecal (-l))\simeq R^1 f_*(\beta^*\Ecal (-lE))\otimes \Ocal_{\pdop^1}(-l)$$
by the projection formula.
We consider the resolution of $E$ 
$$0\ra\Ocal_{\hat Z}(-(j+1)E)\ra\Ocal_{\hat Z}(-jE)\ra(\Ocal_E\otimes\Ocal_{\hat Z}(-E))^{\otimes j} \mbox{\quad (*) \quad}. $$
$\Ocal_E\otimes\Ocal_{\hat Z}(-E)$ is the relative tautological bundle $\Ocal_{\Ncal\dual}(1)$ of $E=\pdop(\Ncal\dual)\ra C$,
with $\Ncal\dual=\Ncal_{C/Z}=L_C^{-1}\oplus L_C^{-1}$, $E\simeq C\times \pdop^1$. 
\\Hence $\Ocal_{\Ncal\dual}(1)=\beta^*L^{-1}\otimes f^*\Ocal_{\pdop^1}(1)$, and tensoring with $\Ecal(j)$ gives
$$\beta^* \Ecal\otimes (\Ocal_E\otimes\Ocal_{\hat Z}(-E))^{\otimes j} =\beta^*\Ecal(-j)_C\otimes f^*\Ocal_{\pdop^1}(j).$$
Hence, after tensoring (*) with $\Ecal(-j)$ and applying $R^i f_*$, we have
\begin{eqnarray} f_*(\beta^* \Ecal\otimes (\Ocal_E\otimes\Ocal_{\hat Z}(-jE))) & = & 0 \mbox{\quad and} \nonumber\\
             R^1 f_*(\beta^* \Ecal\otimes (\Ocal_E\otimes\Ocal_{\hat Z}(-jE))) & = & H^1(\Ecal(-j)_C)\otimes\Ocal_{\pdop^1}(j), \nonumber
\end{eqnarray}
so we have the exact sequence 
$$0\ra R^1 f_*(\beta^* \Ecal(-(j+1)E))\otimes \Ocal_{\pdop^1}(-l)\ra R^1 f_*\beta^*\Ecal(-j E)\otimes\Ocal_{\pdop^1}(-l)\ra \dots $$
$$ \dots \ra H^1(\Ecal(-j)_C)\otimes\Ocal_{\pdop^1}(-(l-j)) \mbox{\quad for \quad}l\le 0.$$
Taking global sections for $j<l$, we finally get 
$$H^0(R^1 f_*(\beta^* \Ecal(-(l+1)E))\otimes (\Ocal_{\pdop^1}(-j)))\simeq H^0(R^1 f_*\beta^*\Ecal(-i E)\otimes\Ocal_{\pdop^1}(-j)),$$
thus especially for $j=l-1$ resp. $j=0$ 
$$H^0(R^1 f_*(\beta^*\Ecal(-l)))\simeq H^0(\Hcal(-l)) \mbox{\quad for \quad}l\le 0,$$
which proves the claim.
\end{prf}

\begin{de}The bundle $\Hcal=\bigoplus_{j=a}^{b}\Ocal(j)^{m_j}$ (with $a\le b, m_a,m_b\neq 0$) on $\pdop^1$ is called the spectrum of 
$\Ecal$. \end{de}

\subsection{Properties of the spectrum}

In the classical case of $\pdop^3$, the spectrum is symmetric (i.e. $\Hcal\dual\simeq \Hcal$ or $\Hcal\dual\simeq \Hcal(-1),$
respectively). Also in our case,  Serre duality induces an isomorphism $\Hcal(\Ecal)\dual\simeq \Hcal(\Ecal\dual(-1))$.
But different normalizations fulfilling 
the condition $1-\mu(L)\le\mu(\Ecal)<0$ are possible, and $\det\Ecal$ is in general not a multiple of $L$. Of course, if there is a 
normalization such that $\det\Ecal\simeq L^n$ for some $n$, the spectrum given by this normalization is symmetric; if 
$\rho=\rk \Pic (Z)=1$, this holds naturally for all possible spectra.
The second classically known property of the spectrum, the connectedness, holds (with small restrictions) also in the case of elliptic 
families:
\begin{thm} 
Assume that for a generic $S\in |\Lambda|$ the linear system $|K_S^{-1}|$ is ample.
\\(i)  \hspace*{0.2cm} If $m_b\neq 0$ for $b>1$, then $m_1,\dots,m_{b-1}\neq 0$.
\\(ii) If $m_a\neq 0$ for $a<-2$, then $m_-2,\dots,m_{a+1}\neq 0$.
\end{thm} 
\begin{prf}
The proof requires three lemmas:
\begin{lem}[Bilinear map lemma]
Let $\phi :U\otimes N\ra N'$ ($U,N,N'$ finite-dimensional vector spaces) be injective on $\{ x\otimes y\}$ for all $x\in U$, then 
$\dim N'-\dim N\ge \dim U -1$.
\end{lem}
\begin{prf} \cite{Ha}, Lemma 5.1 
\end{prf}
\begin{lem}
If $|K_S^{-1}|$ is ample and $ \dim |K_S^{-1}|\ge 2$, then  $|K_S^{-1}|$ contains a pencil of irreducible divisors.
\end{lem}
\begin{flushright}$\Box$\end{flushright}
\begin{lem}
Let $N_{\bullet}\subset H^1(\Ecal_S(\bullet))$ be a graded submodule. Then for $j>0$, either $N_j=0$ or $\dim N_j-\dim N_{-j-1}>0$.
\end{lem}
\begin{prf}
Adjunction formula gives us $K_S^{-1}$. Application of the previous lemma guarantees the existence of $U\subset H^0(L_S)$, such that 
$\dim |U|\ge 1$ and all divisors of $|U|$ are irreducible. 
\\We now consider the multiplication map $U\otimes N_{-j-1}\ra N_j$. There are two possible cases:
\\(a) The map is injective on $\{h\otimes x\}$ for all $h$, hence by the bilinear map lemma $\dim N_j-\dim N_{-j-1}>0$.
\\(b) There exists a $h_0$, such that $N_{-j-1}\rar{h_0} N_j$ is not injective. Since $N_{-j-1}\rar{h} N_j$ is injective for generic $h$,
we have $\dim N_j\ge\dim N_{-j-1}$. If we denote $C=\di (h_0)$, we have the following diagram:
$$\begin{diagram} \node {0} \arrow{e} \node{H^0(\Ecal_C(-j))} \arrow{e} \node{H^1(\Ecal_S(-j-1))} \arrow{e,t}{h_0} \node{H^1(\Ecal_S(-j))}
\\ \node{0} \arrow{e} \node{M_j} \arrow{n,l}{\subset} \arrow{e} \node{N_{-j-1}} \arrow{n,l}{\subset} \arrow{e,t}{h_0} \node{N_j}
\arrow{n,l}{\subset}
\end{diagram}$$
with 
$$M_{\bullet}=\ker(H^1(\Ecal_S(\bullet))\rar{h_0} H^1(\Ecal(\bullet))).$$
Hence $\dim(N_{-j-1}/M_{-j})\le\dim (N_{-j}/M_{-j+1})$, therefore $\dim N_{-j} - \dim N_{-j-1} \ge \dim M_{-j+1} - \dim M_{-j}$.
But $C$ is irreducible, hence $M\subset H^0(\Ecal_C(\bullet))$, so we can again apply the bilinear map lemma to conclude 
$\dim M_{-j+1} - \dim M_{-j}>0$, which finishes the proof.
\end{prf}
Because of $\Hcal(\Ecal)\dual\simeq \Hcal(\Ecal\dual(-1))$ it is sufficient to proof $a)$.
We consider the sequence \Shoexaseq{\Ecal(-l-1)}{\Ecal(-l)}{\Ecal_S(-l)} and apply the previous lemma to the long exact sequence,
beginning with 
$$0\ra H^1(\Ecal{-l-1}) \ra H^1(\Ecal(-l)) \onto N_{-j}\subset H^1(\Ecal_S(-l)).$$
Evaluating this via the spectrum, we have in the left term 
$$h^1(\Ecal{-l-1})=h^0(\Hcal(-l-1))=\sum m_k h^0(\Ocal_{\pdop^1}(k-l-1))=\sum_{k\ge l+1}m_k(k-l)$$
and in the middle term
$$h^1(\Ecal{-l})=h^0(\Hcal(-l))=\sum m_k h^0(\Ocal_{\pdop^1}(k-l))=\sum_{k\ge l}m_k(k-l+1),$$
hence $\dim N_{-l}=\sum_{k\ge j} m_k$ and $\dim N_{-l} - \dim N_{-l-1}=m_j$. Application of lemma 4.7 finishes the proof.
\end{prf}
The connectedness allows us to restrict the maximal und minimal summands in the $\pdop^1$-Bundle $\Hcal$ and to proof a vanishing
theorem:
\begin{cor}
With $d=\deg \Hcal$, $r=\rk \Hcal$ and $\Hcal=\bigoplus_{j=a}^{b}\Ocal(j)^{m_j}$, the inequalities 
$$a\ge \frac{d-1}{r+2} - \frac{r+1}{2} \mbox{\quad and \quad} b\le\frac{d-1}{r+2} + \frac{r+1}{2}$$
hold.
\end{cor}
\begin{cor} $H^1(\Ecal(-l))=0$ for $l\ge\min(0,\frac{d-1}{r+2}+\frac{r+1}{2})$, 
$H^2(\Ecal(l))=0$ for $l\ge\max(-1,\frac{d-1}{r+2}-\frac{r+3}{2})$.
\end{cor}
\begin{prf} The second corollary is an immediate conclusion from the first. To prove the first,
we consider the function $f(l)=\sum_{\nu\le l}m_{\nu}$. It is obvious from the definition, that
\\(i)  $0\le f(l) \le r$
\\(ii) $f(l)=0$ for $l<a$
\\(iii)$f(l+1)>f(l)$ for $a-1\le l \le b$ and $l\neq -1,-2$.
\\Now, we have 
\begin{eqnarray} d & = & \sum_{l\le N} l f(l)-f(l-1) \mbox{\quad for \quad} N>> 0 \nonumber \\
                   & = & \sum_{l\le N} l f(l) - \sum_{l<N} (l+1)f(l) \nonumber \\
                   & = & Nr - \sum_{l<N} f(l) \nonumber \\
                   & = & Nr - \sum_{l\le-2} f(l)-f(-1)-f(0)-\sum_{0<l<N} f(l) \nonumber 
\end{eqnarray}
We assume $-r<a\le -2$ (otherwise, the inequalities are trivially true) and split of the equation to 
apply (iii):
\begin{eqnarray} d & = & Nr-\sum_{a\le l\le -2} (f(l)-(l-a+1))-\sum_{0<l<a+r} (f(l)-(l-a-1)) - \nonumber \\
                   &   &-f(0)-f(-1)-\sum_{\nu=1}^{r-1} \nu -\sum_{l=a+r+1}^{N-1} f(l) \nonumber \\
                   & \le & Nr -r(N-a-r-1)+2(a+1) -\frac{r(r-1)}{2} \nonumber
\end{eqnarray}
The proof for the second inequality is analogous.
\end{prf}

\end{document}